\documentclass[11pt]{article}

\usepackage[margin=1.5in]{geometry}
\usepackage{amsmath, amssymb,amsfonts}

\usepackage{listings}
\lstdefinelanguage{mathlang} {
  literate=
    {→}{{\ensuremath{\rightarrow}}}1
    {θ}{{\ensuremath{\mathrm{\theta}}}}1
    {∀}{{\ensuremath{\forall}}}1
    {∃}{{\ensuremath{\exists}}}1
    {λ}{{\ensuremath{\mathrm{\lambda}}}}1
    {₀}{{\ensuremath{_0}}}1
    {₁}{{\ensuremath{_1}}}1
    {₂}{{\ensuremath{_2}}}1
    {ℕ}{{\ensuremath{\mathbb{N}}}}1
    {∧}{{\ensuremath{\wedge}}}1
    {≡}{{\ensuremath{\equiv}}}1
    {⟶}{{\ensuremath{\longrightarrow}}}1
    {≤}{{\ensuremath{\leq}}}1
    {∈}{{\ensuremath{\in}}}1,
  basicstyle=\ttfamily\small
}
\newcommand{\ZZ}{\mathbb{Z}}
\newcommand{\CC}{\mathbb{C}}
\newcommand{\ph}{\varphi}
\newcommand{\ttt}[1]{\texttt{#1}}

\title{Mathematics and Language}

\author{Jeremy Avigad}

\begin{document}

\maketitle

\begin{quote}

I learned empirically that this came out this time, that it usually does come out; but does the proposition of mathematics say that? \ldots\ The mathematical proposition has the dignity of a rule.

\emph{So} much is true when it's said that mathematics is logic: its moves are from rules of our language to other rules of our language. And this gives it its peculiar solidity, its unassailable position, set apart.

--- Ludwig Wittgenstein

\bigskip

\ldots\ it seemed to me one of the most important tasks of philosophers to investigate the various possible language forms and discover their characteristic properties. While working on problems of this kind, I gradually realized that such an investigation, if it is to go beyond common-sense generalities and to aim at more exact results, must be applied to artificially constructed symbolic languages.\ldots\ Only after a thorough investigation of the various language forms has been carried through, can a well-founded choice of one of these languages be made, be it as the total language of science or as a partial language for specific purposes.

--- Rudolf Carnap

\bigskip

Physical objects, small and large, are not the only posits.\ldots\ the abstract entities which are the substance of mathematics\ldots\ are another posit in the same spirit. Epistemologically these are myths on the same footing with physical objects and gods, neither better nor worse except for differences in the degree to which they expedite our dealings with sense experiences.

--- W.~V.~O.~Quine

\bigskip

``When I use a word,'' Humpty Dumpty said in rather a scornful tone, ``it means just what I choose it to mean --- neither more nor less.''

``The question is,'' said Alice, ``whether you can make words mean so many different things.''

``The question is,'' said Humpty Dumpty, ``which is to be master --- that's all.''

--- Lewis Carroll

\end{quote}

\section{Introduction}

Mathematics holds a special status among the sciences. We expect the results of mathematical calculation to tell us things about the world, such as where to look for Mars or Venus at a given date and time, whether a bridge is strong enough to withstand the loads that will be placed on it, and whether a straight is more likely than a flush in poker. At the same time, we do not look to the world to confirm our mathematical knowledge: we don't rely on empirical observations to justify the claim that two and two make four, and it is hard to imagine a scientific experiment that could disconfirm that belief. 

There are other features that set mathematics apart. Unlike physical objects, mathematical objects are not located in space. It is almost grammatically incorrect to speak about mathematical objects in temporal terms, for example, to say that seven was prime yesterday or that it will, in all likelihood, be prime tomorrow. And it seems that, at least ideally, mathematical knowledge admits a kind of certainty that cannot be attained empirically. Our mathematical claims are sometimes mistaken, but that is generally attributed to a failure to understand the proper methods of reasoning or follow them correctly. In contrast, the doubts we have about empirical claims seem to be inherent in the nature of empirical methods, and not in our ability to carry them out.

The philosophy of mathematics tries to explain these features: the \emph{ontology of mathematics} aims to provide an account of mathematical objects and their nature, and the \emph{epistemology of mathematics} aims to provide an account of the means by which we can reliably come to have mathematical knowledge. Whatever we have to say about the nature of mathematical objects, an important part of the game is to explain how it is we can possibly know things about them. And whatever we have to say about that knowledge, an important part of the game is to explain why it is knowledge worth having, and, in particular, how it grounds rational behavior and expectations with respect to empirical phenomena.

Philosophers have played this game since the origins of philosophy itself. Plato relied on his theory of forms, while Aristotle invoked a rival, less Platonic, theory of forms. Descartes' \emph{Meditations} explains why we can be nearly certain that our mathematical claims are true and have something to do with the world, though they are not quite as certain as the fact that God exists. Leibniz classified mathematical truths as analytic truths, Locke classified them as ideas of relation, and Hume classified them as relations of ideas. For Kant, arithmetic and geometric truths were important examples of knowledge that is synthetic (roughly, not ``contained'' in the definitions of the constituent concepts) and \emph{a priori} (epistemologically independent of our experience of the world). His momentous \emph{Critique of Pure Reason} sought to explain how such knowledge is possible.

Without getting bogged down in knotty metaphysics, modern logic gives us \emph{descriptive} accounts of how mathematics works. Various formal axiomatic foundations --- set-theoretic foundations, type-theoretic foundations, categorical foundations --- provide us with models of the mathematical language and its general principles of inference. To be sure, there are puzzling features of mathematical language; for example, we sometimes refer to ``the seven-element cyclic group'' knowing full well that it is only defined up to isomorphism. Different frameworks account for such features in different ways, but they are generally inter-interpretable, and taken together, they give us a pretty good low-level account of the structure of mathematical language and inference.

There are also fairly lowbrow explanations as to why mathematics exhibits the features enumerated above. Mathematical objects are what the norms of mathematical practice say they are, no more, and no less, and we come to have mathematical knowledge by following the norms of reasoning correctly, not by performing experiments and measurements. We learn mathematics by learning the appropriate norms from our parents, our teachers, and the mathematical literature, with additional feedback and correction from graders and journal referees. 

But why are the norms of mathematical practice what they are, and how is it that following those norms tells us anything at all about the way the world is? In the first half of the twentieth century, the logical positivists viewed these questions as two sides of the same coin. Mathematics tells us something about the world precisely because it embodies norms and conventions with which we describe the world and make sense of physical experience. And we have settled on these norms just because they provide us with useful ways of making sense of our experiences; if they weren't useful, we would adopt different norms. 

On this view, what makes something a piece of mathematics proper is that we hold the rules and conventions relatively fixed in our ordinary patterns of reasoning, blaming any mismatches with experience on the application rather than the rules. In a sense, there are plenty of counterexamples to the mathematical claim that one and one make two: one drop of water and another drop of water make a bigger drop of water, one rabbit and another rabbit make lots of rabbits, and one rock combined with another rock with sufficient force makes dust. But in each case, we reject the counterexample as a misapplication of the mathematical principle. We hold our norms for counting and computing sums fixed because they provide generally useful ways of interpreting our experiences, but in ordinary circumstances we keep experience beholden to the rules, and not the other way around. This point of view has been nicely summed up by William Tait:
\begin{quote}
Only empirical explanation is possible for why we have come to accept the basic principles that we do and why we apply them as we do---for why we have mathematics and why it is at it is. But it is only within the framework of mathematics as determined by this practice that we can speak of mathematical necessity. In this sense, which I believe Wittgenstein was first to fully grasp, mathematical necessity rides on the back of empirical contingency.
\end{quote}

For the philosopher Rudolf Carnap, adopting a set of mathematical norms and conventions was part of choosing a suitable ``language form,'' a task which he believed could be informed by philosophical reflection and guidance. Carnap and his fellow logical positivists drew a distinction between ``analytic'' truths, like mathematical truths, that are grounded in the norms of the linguistic framework, and ``synthetic'' truths, which are grounded by empirical observations. In his influential essay, ``Two dogmas of Empiricism,'' W.~V.~O.~Quine famously denied that there is a sharp distinction between the conventional and empirical parts of language: any empirical claim can be held true in the face of recalcitrant experience by modifying other parts of our language accordingly, and, conversely, even our most deeply held mathematical beliefs are subject to revision in light of future experiences. Thus, for Quine, \emph{all} of science is a matter of maintaining a web of beliefs and making it fit with our experiences. On this view, the allegedly special character of mathematical claims as unassailably fixed is an illusion; at best, we can say that mathematical claims are more stable, or less subject to revision, than their ostensibly empirical counterparts.

In support of Quine's criticisms, one can observe that mathematical concepts do indeed evolve over time, though the subject seems to be remarkably adept at reinterpreting prior developments in a way that preserves their essential correctness. In any case, whether we view the outward features of mathematical practice as reflective of a special status of its claims or simply a provisional stance we adopt towards them, there remains a substantial point of agreement in the views of Carnap and Quine: mathematical and scientific reasoning are shaped by the rules of our language, and these rules are, in turn, adopted for pragmatic scientific reasons.

My goal in this essay is to take this point of view seriously, and argue that it is an important task for philosophy to clarify the rules and norms of our mathematical practices and to help us understand the role they play in our reasoning. For all their talk about the general considerations that influence our choice of language, neither Carnap nor Quine was very specific about how to assess the effects of the choices that we have made historically and continue to make in our everyday activities. Note that I am using the word ``language'' here in a very broad sense to refer not only to our choice of words and grammar, but also the way we frame, conceptualize, analyze, and reason about the problems before us.

What I am advocating is a view of mathematics as a linguistic artifact, something we have designed, and continue to design, to help us get by in the world. As a community of practitioners, we choose to do mathematics the way we do, and when we do mathematics in the usual ways, we are bound by these choices. Philosophy can help us become aware of the choices we have made, and understand why we have made them, and even, perhaps, why we might want to change them. Various pieces of mathematics are designed to help us achieve various goals, and like any human artifacts, they can serve their purpose well or poorly. The philosophical challenge is to understand the general principles by which they serve our goals well, and to support the development of better artifacts. In that sense, the philosophy of mathematics can be viewed as a kind of language engineering, akin to any other design science.

In this essay, I will try to pursuade you that this view of mathematics makes sense, and provides us with a framework for understanding aspects of mathematics that are important to us. In the next section, I will share some reflections on ordinary mathematical language, using the history of mathematics as a guide. I will then turn to formal mathematical languages, which is to say, symbolic languages whose grammar and rules of use are expressed in precise mathematical terms. Finally, I will explore some of the relationships between the two, and use those observations to bolster my recommendations for the philosophy of mathematics. 

\section{Informal mathematical language}
\label{section:informal}

I will begin with some observations regarding ordinary mathematical concepts and the way they evolve over time. In 2010, a graduate student in my department, Rebecca Morris, asked me to supervise her research on the history and philosophy of mathematics. She was particularly interested in the development of the concept of a ``function,'' and hoped that studying the history would help illuminate the way we reason about functions in everyday mathematics. From a philosophical standpoint, she wanted to understand \emph{why} we talk about functions the way we do, and the senses in which we are justified in doing so.

Today, we think of a function as a correspondence between any two mathematical domains: we can talk about functions from the real numbers to the real numbers, functions that take integers to integers, functions that map each point of the Euclidean plane to an element of a certain group, and, indeed, functions between any two sets. As we began to explore the topic, Morris and I learned that most of the historical literature on the function concept focuses on functions from the real or complex numbers to the real or complex numbers. There is a good reason for that: for most of the history of the function concept, from Euler's 1748 introduction to the analysis of the infinite to the dramatically different approaches to complex function theory by Riemann and Weierstrass in the latter half of the nineteenth century, the word ``function'' was used exclusively in that sense. Moreover, authors generally tended to assume, implicitly, than any function is given by an analytic expression of some sort. In 1829, in a paper on the Fourier analysis of ``arbitrary functions,'' Dirichlet went out of his way to reason about functions extensionally, which is to say, without making use of any such representation. In that paper, he famously put forth the example of a function that takes one value on rational numbers and another value on irrational numbers as an extreme case, one that was difficult to analyze in conventional terms.

Dirichlet's ``arbitrary'' functions were, nonetheless, functions from the real numbers to the real numbers. Even the notion of a ``number theoretic function,'' like the factorial function or the Euler function, is nowhere to be found in the literature; authors from Euler to Gauss referred to such entities as ``symbols,'' ``characters,'' or ``notations.'' Morris and I tracked down what may well be the first use of the term ``number theoretic function'' in a paper by Eisenstein from 1850, which begins with a lengthy explanation as to why it is appropriate to call the Euler phi function a ``function.'' We struggled to parse the old-fashioned German, which translates roughly as follows:
\begin{quote}
Once, with the concept of a function, one moved away from the necessity of having an analytic construction and began to take its essence to be a tabular collection of values associated to the values of one or several variables, it became possible to take the concept to include functions which, due to conditions of an arithmetic nature, have a determinate sense only when the variables occurring in them have integral values, or only for certain value-combinations arising from the natural number series. For intermediate values, such functions remain indeterminate and arbitrary, or without any meaning. 
\end{quote}
When the gist of the passage sank in, we laughed out loud. The opening phrase is a nod to Dirichlet's notion of a function as an arbitrary tabular pairing of input and output values. We had expected Eisenstein to go on to say something like this: ``Of course, there is no reason to restrict attention to the real and complex numbers; we can now take a function to be a tabular pairing of values in any two domains.'' Instead, he went on to argue that we can extend the function concept to functions defined on the integers by considering them to be \emph{partially defined functions on the real numbers}, that is, functions from the real numbers to the real numbers that happen to take values only on integer arguments. Indeed, the first hints of the modern notion of a function as an arbitrary correspondence between any two domains did not appear until the late 1870's, and Eisenstein's introductory paragraph is a nice reminder of the fact that mathematical concepts evolve gradually, and often in surprising ways.

Morris and I decided to focus our attention on another type of function that was studied in the nineteenth century (though they were not labeled as such). In 1837, Dirichlet proved a beautiful theorem that states that there are infinitely many prime numbers in any arithmetic progression in which the first two terms do not have a common factor. For example, there is clearly only one prime number in the progression $2, 12, 22, \ldots$ because the first two terms, and hence all the terms, are multiples of $2$. Dirichlet's theorem states that this is essentially the only thing that can go wrong. For example, in the progression $3, 13, 23, 33, \ldots$, the first two terms have no factor in common, and so there will be infinitely many primes. This was conjectured to be true around the turn of the nineteenth century but even the great master, Gauss himself, had been unable to prove it. Dirichlet's proof is a landmark not only because it solved a hard problem, but also because it introduced methods that are now central to algebraic and analytic number theory.

Modern presentations of Dirichlet's proof rely on the notion of a \emph{Dirichlet character}. To deal with the arithmetic progression with first term $m$ and common difference $k$, we consider the multiplicative group $(\ZZ/k\ZZ)^*$ of invertible residues modulo $k$, that is, the numbers modulo $k$ that have no common factor with $k$. For example, when $k = 10$, we consider the set of numbers $\{1, 3, 7, 9\}$ with multiplication modulo $10$. A \emph{character} on this group is a nonzero homomorphism from this group to the complex numbers, that is, a function $\chi$ with the property that $\chi(x \cdot y) = \chi(x) \cdot \chi(y)$ for every $x$ and $y$. For example, there is the trivial character on $(\ZZ/10\ZZ)^*$ which sends all of the residues to $1$, and another character $\chi$ which sends $1, 3, 7, 9$ to $1, i, -i, -1$, respectively. There are four distinct characters modulo $10$, and one can show that, for any $k$, there are only finitely many characters modulo $k$. To be precise, a \emph{Dirichlet character} is obtained from an ordinary (group) character by ``lifting'' the values to the integers, i.e.~defining $\hat \chi(n)$ to be the value of $\chi$ at $n$ modulo $k$ if $n$ has no common factor with $k$, and $0$ otherwise. As is common practice, I will blur the distinction between the group character and the associated Dirichlet character, and use the term ``character'' for both.

What made Dirichlet's theorem interesting to us is that in contemporary proofs, characters are treated as ordinary mathematical objects like the natural numbers. For example, in an ordinary textbook proof:
\begin{enumerate}
\item we show that the characters modulo $k$ form a group under pointwise multiplication, just as the integers modulo $k$ form a group under addition;
\item we define functions $L(s, \chi)$, called \emph{Dirichlet L-series}, whose second argument is a character;
\item we write down sums $\sum_\chi \chi(m)^{-1} L(s, \chi)$ where the index $\chi$ ranges over characters, just as we write down sums like $\sum_{i = 1}^{n} i^2$ where the index $i$ ranges over the integers.
\end{enumerate}
These features are distinctly modern: mathematicians from Euler to Cauchy seemed to think of functions as syntactic expressions, like $f(x) = x^2$, that are somewhat different from objects like numbers. In today's terms, functions are ``higher-order objects,'' and what is notable is that in contemporary mathematics their status is no different from more fundamental mathematical entities. Morris and I reasoned that if we could understand how the features above were implemented in Dirichlet's original proof, it would help us understand how functions attained their contemporary status. 

When we read Dirichlet's original papers, however, we were surprised to discover that there is no notion of character there at all! Rather, there are certain algebraic expressions which, today, we view as values of the characters. Where we would write $\chi(n)$, Dirichlet wrote expressions
\[
\theta^\alpha \ph^\beta \omega^\gamma \omega'^{\gamma'} \ldots,
\]
where $\theta, \ph, \omega, \omega', \ldots$ are complex roots of $1$ and $\alpha, \beta, \gamma, \gamma', \ldots$ are certain integers that implicitly depend on $n$. Roughly, whereas today we define the characters as the nonzero homomorphisms from the group of invertible residues to the complex numbers, Dirichlet used an explicit symbolic representation of the characters and implicitly relied on the fact that they have that key property.

Dirichlet showed that the relevant symbolic expressions can be parameterized by tuples $\mathfrak{a}, \mathfrak{b}, \mathfrak{c}, \mathfrak{c}', \ldots$ of integers, and wrote $L_{\mathfrak{a}, \mathfrak{b}, \mathfrak{c}, \mathfrak{c}',\ldots}(s)$ where we would write $L(s,\chi)$. Moreover, where we would write a summation over characters, Dirichlet summed over the representing tuples. For example, where we would write
\[
\sum_{\chi \in \widehat{(\mathbb{Z} / k \mathbb{Z})^*}}\overline{\chi(m)}\log L(s, \chi),
\]
Dirichlet wrote 
\[
\sum \Theta^{-\alpha_{m}\mathfrak{a}}\ \Phi^{-\beta_{m}\mathfrak{b}}\Omega^{-\gamma_{m}\mathfrak{c}}\Omega^{-\gamma'_{m}\mathfrak{c}'} \cdots \log L_{\mathfrak{a},\mathfrak{b},\mathfrak{c},\mathfrak{c}', \ldots.}
\]

Dirichlet's reliance on representing data played out in other ways. In contemporary presentations of Dirichlet's proof, we divide the characters into three groups: first, there is the trivial character, $\chi(m) = 1$ for every residue $m$; then there are the characters, other than the trivial character, that take only the real values $\pm 1$; and, finally, the remaining characters, each of which takes on at least one non-real complex value. This categorization is \emph{extensional}, which is to say, it is cast in terms of the \emph{values} of the characters $\chi$, and not the data used to represent them. In contrast, Dirichlet's division into cases was cast in terms of the defining data, something we would describe as an \emph{intensional} categorization.

While treating characters as first-class mathematical objects, modern proofs also make them objects of study in their own right. Morris and I were curious to understand how the understanding of Dirichlet characters evolved over time, and so we studied a number of subsequent presentations of Dirichlet's theorem and related results: extensions, by Dirichlet, of his result to Gaussian integers and quadratic forms; an 1863 presentation by Dedekind, in a supplementary appendix to his writeup of Dirichlet's lecture notes on number theory; generalizations by Dedekind, in 1879, and Weber, in 1882, of the notion of a character to an arbitrary group;  a constructive proof, due to Kronecker, presented in the 1880's and later written up by his student, Hensel; presentations of Dirichlet's theorem and extensions by Hadamard in 1896 and de la Vall\'ee Poussin in 1897; and presentations by Landau in 1909 and 1927, the latter which reads much like a contemporary textbook proof.

What we found is that the path from Dirichlet's proof to Landau's was long and tortuous, with various ``modern'' features appearing hesitatingly and tentatively, in fits and starts. Over time, authors gradually identified the expressions defining the characters as key components in the proofs, separated out their central properties, and proved them as independent lemmas. One advantage to doing this is that the fiddly technical details of manipulating representations become isolated and black-boxed, so that readers only have to remember the key properties, and not their proofs. Another advantage is that identifying the notion of a character and their essential properties supports the generalization of the notion to other groups, as well as generalizations of Dirichlet's proof itself. Dirichlet himself made moves in this direction in his later work. Dedekind's 1863 presentation used $\chi$ to denote the values of the characters and separated out some of their axiomatic properties, and in 1882, Weber gave the general definition of a character of an abelian group, and proved the general properties. In 1909, Landau enumerated four ``key properties'' of the characters, and emphasized that these are all that is needed in the proof of Dirichlet's theorem; once we have them, we can forget all the details of the representations of the characters.

The treatment of characters on par with common mathematical objects like numbers came gradually as well. Where we would write sums over the characters, Dedekind's 1863 presentation summed over representing data with expressions like $\sum_{a, b, c, c', \ldots}$, like Dirichlet's. In 1897, Hadamard followed a different approach: he introduced an arbitrary numbering of the characters, writing them $\psi_1, \psi_2, \ldots, \psi_M$, which meant that he could sum over the characters with conventional notation $\sum_{i = 1}^M$. In contrast, de la Vall\'ee Poussin introduced special notation $S_\chi$ to sum over characters, while at the same time using the notation $\sum_i$ for summing over a finite set of integers. This suggests an abstraction of the summation operation to the set of characters, but marked by a concern that summing over a finite set of characters is nonetheless an operation that is distinct from summing over a finite set of integers. In 1909, Landau followed Hadamard's approach, but in 1927 he finally wrote $\sum_\chi$, as we do today.

The transition to treating the characters as arguments to $L$-functions was similarly prolonged. For example, Hadamard simply wrote $L_v$ for the series corresponding to $\psi_v$. In 1897, de la Vall\'ee Poussin wrote $Z(s, \chi)$, and Landau adopted the notation $L(s, \chi)$ that we use today. And the way authors described the classification of $L$-series, which plays a crucial role in the proof, also changed over time. Most authors we considered favored an intensional characterization like Dirichlet's, though some authors gave both characterizations, that is, specifying the division in terms of the output values as well as in terms of the representing data.

In short, it took almost 90 years for the treatment of Dirichlet characters to evolve to the point where they were identified as objects of interest, seen to have a group structure, passed as arguments to functions, taken to be the indices of a summation, and so on. Since then, from Landau's 1927 proof to the present day, the treatment of Dirichlet characters has remained remarkably stable.

The point I wish to make is that these changes constitute an ontological shift, pure and simple. In 1837, Dirichlet could not say ``a character is a nonzero homomorphism from $(\ZZ/\ZZ_k)^*$ to $\CC$,'' because there was no general notion of a function between two domains. He certainly could not do all the things that Landau could do, because the linguistic and methodological resources were simply not available to him. Here I am once again using the term ``language'' in a broad sense, to include not just grammar and vocabulary but also the conceptual and inferential apparatus that determine the way we talk about and reason about a subject.

It should not be surprising that such changes occur only gradually. Whenever we write down a piece of mathematics, it is with the intention that others will read it, understand it, and judge it to be correct and worthwhile. When we introduce new terminology, notation, and patterns of inference, it is usually the case that the intended usage can be explained and justified in conventional terms. But when the changes represent a marked departure from the status quo, they bring with them a host of nagging concerns. What do the new terms and symbols mean? What are the rules governing their use? Are these rules reasonable, and justified? Are they appropriate to the subject matter, and do they give us legitimate answers to the questions we have posed? It is not sufficient that a single author becomes comfortable with the changes; for mathematics to work, the entire mathematical community has to function as a coherent body and come to agreement as to what is allowed. Thus, when changes occur, they occur for good reasons. The ultimate goal of mathematics is to get at the truth, and when new methods make it possible to push the boundaries of our knowledge, there is great incentive to make sense of them, get used to them, and incorporate them into the canon.

Our study of the history of Dirichlet's theorem showed that there are, indeed, good reasons to treat characters as ordinary mathematical objects. Simply identifying characters as objects of interest makes it possible to highlight the properties that make then useful. The very act of naming them simplifies the expressions that occur in complex proofs, and developing their properties independently makes those proofs more modular. Dependencies between pieces of a mathematical text are minimized, so that one can verify properties of the characters independently, and then suppress those details later on. This makes it easier to read the proof, understand it, and verify its correctness.

Such a modular structuring has additional benefits, in that the modularized components are adaptable, reusable, and generalizable. Having a theory of characters at one's disposal means having a tool that can be used to solve other problems, and having an abstract description of what makes Dirichlet characters useful enables us to design or discover other objects with similar properties. Recall that Dedekind and Weber transferred the notion of character to arbitrary finite groups. Later generalizations to continuous groups and nonabelian groups are now fundamental to the subjects of representation theory and harmonic analysis.

Moreover, treating characters on par with other mathematical objects allows other bodies of mathematics to be invoked and applied uniformly. Mathematicians were already comfortable summing over finite sets of numbers; if characters have the same status as numbers, conventional notation and methods of calculation transfer wholesale. Similarly, by the 1870's, the mathematical community began to gain facility with the abstract notion of a group. If one can consider a group of characters in the same manner as a group of residues, both are subsumed under a common body of knowledge.

The history of Dirichlet's theorem is only one manifestation of the momentous changes that swept across mathematics in the nineteenth century. The development of the modern theory of functions and the introduction of ideals and quotients in algebra are additional forces that pushed for both algebraic and set-theoretic abstraction --- a new ``language form,'' in Carnap's terminology. These are the sorts of changes that we need to study carefully, because understanding them is crucial to understanding how we do mathematics today. In the next two sections, I will argue that we can be more scientific in this study, and that formal languages and logical methods can play an important role.

\section{Formal mathematical languages}
\label{section:formal}

In 1931, Kurt G\"odel, then still a relatively unknown young logician, published his two \emph{incompleteness theorems}. These are now recognized as among the most profound and important results in the history of logic. The first incompleteness theorem shows, roughly speaking, that no formal axiomatic system of mathematics can settle all mathematical questions, even questions about the natural numbers; there will always be relatively straightforward statements about the natural numbers that cannot be proved or disproved on the basis of the axioms. The second incompleteness theorem shows that, moreover, no such system can establish its own consistency. The opening sentences of G\"odel's paper are compelling even today, delivering the dramatic news in a calm, matter-of-fact tone:

\begin{quote}
The development of mathematics toward greater precision has led, as is well known, to the formalization of large tracts of it, so that one can prove any theorem using nothing but a few mechanical rules. The most comprehensive formal systems that have been set up hitherto are the system of \emph{Principia Mathematica (PM)} on the one hand and the Zermelo-Fraenkel axiom system of set theory (further developed by J.~von Neumann) on the other. These two systems are so comprehensive that in them all methods of proof used today in mathematics are formalized, that is, reduced to a few axioms and rules of inference. One might therefore conjecture that these axioms and rules of inference are sufficient to decide \emph{any} mathematical question that can at all be formally expressed in these systems. It will be shown below that this is not the case\ldots
\end{quote}

What I would like to focus on here are not G\"odel's negative, limitative results, but the positive side of his assessment. Ernst Zermelo presented his system of axiomatic set theory in 1908, and Bertrand Russell and Alfred North Whitehead presented a system of ``ramified type theory'' in their monumental work, \emph{Principia Mathematica}, the three volumes of which first appeared between 1911 and 1914. It is striking that as early as 1931, G\"odel could convey the conventional understanding that, indeed, these systems suffice to formalize the vast majority of ordinary mathematical reasoning, a fact that remains true to this day.

But how do we get to a proof of Fermat's last theorem from a few simple axioms and rules? The phenomenon is analogous to the way that we obtain complex computational systems from primitive components. Conceptually, computer microprocessors and memory chips are built from very simple elements, components that can store a bit of memory or carry out a simple logical operation. These are combined to form components that are capable of carrying out more complex operations, and those operations, in turn, form the basis for low-level programming languages. Working our way up in a modular fashion, we get operating systems and compilers, and ultimately systems that manage airline schedules, control industrial power plants, play chess, and compete in game shows.

The ``evidence'' G\"odel had in mind for the claim that mathematics can be formalized is similar. Starting with a few axioms and rules, we can define basic mathematical objects like numbers, functions, sequences, and so on, and derive their properties. From there, we work our way up through the mathematical canon. In practice, however, the process is rather tedious: we think about mathematics at a higher level of abstraction, and even with a modular approach, there are far too many details to write down in practice. Famously, in the three volumes of formal mathematics in the \emph{Principia}, Russell and Whitehead never got past elementary arithmetic. Explaining how complex objects and inferences can be built up from more simple ones and carrying out a number of examples is enough to make a compelling case that mathematics can be formalized \emph{in principle}, but that is not the same as saying that it can be formalized \emph{in practice}.

Towards the end of the twentieth century, however, computer scientists began to develop ``computational proof assistants'' that now make it possible to construct axiomatic proofs of substantial theorems, in full detail. The languages with which one interacts which such systems are much like higher-level programming languages, except that the ``programs'' are now ``proofs'': they are directives, or pieces of code, that enable the computer to construct a formal axiomatic proof. Proof assistants keep track of definitions and previously proved theorems, as well as notation and abbreviations. They also provide various types of automation to help fill in details or supply information that has been left implicit. Users can, moreover, declare bits of knowledge and expertise as they develop their mathematical libraries, such as rules for simplifying expressions and inference patterns to be used by the automation procedures.

There are a number of proof assistants now in existence, as well as fully verified libraries of number theory, real and complex analysis, measure theory and measure-theoretic probability, linear algebra, abstract algebra, and category theory. Substantial mathematical theorems have been verified in such systems. Proof assistants are also used to verify claims as to the correctness of hardware and software with respect to their specifications, something that is especially valuable when human lives are at stake, as is the case, for example, with control systems for cars and airplanes. Describing the state of the field today would take me too far afield, but there are suggestions for further reading in the notes at the end of this essay.

Designing a theorem prover thus involves designing a language that can express mathematical definitions, theorems, and proofs, providing enough information for the computer to construct the underlying formal mathematical objects, in a manner that is convenient, efficient, and practically manageable. To start with, one has to choose a formal system that constitutes the underlying standard of correctness. Axiomatic set theory, in which every mathematical object is encoded as a set, has proved to be a remarkably flexible and robust foundation, but experience shows that additional structure is useful when it comes to interactive theorem proving. The fact that any expression in set theory can be interpreted as a natural number, a function, or a set makes it hard for the system to infer the user's intent, thus requiring users to provide much more information. In various systems of \emph{type theory}, every object is tagged as an element of an associated \emph{type}. A given expression may denote a natural number, a function from natural numbers to natural numbers, a pair of numbers, or a list of numbers; the syntactic form of the expression, and the manner it is presented, indicates which is the case. Such languages thus provide a closer approximation to ordinary mathematical vernacular.

Mathematical definitions and theorems are then expressed in the chosen foundational language. Here is a statement of the prime number theorem in the Isabelle theorem prover:
\begin{quote}\begin{lstlisting}
theorem PrimeNumberTheorem: "(λn. pi n * ln n / n) ⟶ 1" 
\end{lstlisting}\end{quote}
Here \texttt{pi} denotes the function $\pi(n)$ that counts the number of primes less than or equal to $n$:
\begin{quote}\begin{lstlisting}
pi n ≡ card {p. p ≤ n ∧ p ∈ prime} 
\end{lstlisting}\end{quote}
The Feit-Thompson theorem is a landmark theorem of finite group theory; it asserts that every finite group of odd order is solvable, or, equivalently, that the only simple groups of odd order are cyclic. These two assertions are expressed as follows in the SSReflect dialect of Coq:
\begin{quote}\begin{lstlisting}
Theorem Feit_Thompson (gT : finGroupType) (G : {group gT}) :
  odd #|G| → solvable G.

Theorem simple_odd_group_prime (gT : finGroupType) 
    (G : {group gT}) :
  odd #|G| → simple G → prime #|G|. 
\end{lstlisting}\end{quote}
The Kepler conjecture asserts that there is no way of packing equally-sized spheres into space that results in a higher-density than that achieved by the familiar means of arranging them into nested hexagonal layers, the way oranges are commonly stacked in a crate. Here is a statement of this theorem in HOL Light:
\begin{quote}\begin{lstlisting}
∀V. packing V => (∃c. ∀r. &1 <= r => 
  &(CARD(V INTER ball(vec 0,r))) <=
    pi * r pow 3 / sqrt(&18) + c * r pow 2))
\end{lstlisting}\end{quote}
And here is a statement of the Blakers-Massey theorem, a theorem of algebraic topology, formulated in a formal framework known as \emph{homotopy type theory} and rendered in the Agda proof assistant:
\begin{quote}\begin{lstlisting}
blakers-massey : ∀ {x₀} {y₀} (r : left x₀ ≡ right y₀) → 
  is-connected (n +2+ m) (hfiber glue r)
\end{lstlisting}\end{quote}

All of the theorems just listed have been fully verified in the systems indicated. To that end, all of the objects they mention, such as natural numbers, real numbers, groups, and homotopy fibers, had to be defined in the relevant foundational frameworks. In the Lean theorem prover, for example, we can define the natural numbers as follows:
\begin{quote}\begin{lstlisting}
inductive nat : Type :=
| zero : nat
| succ : nat → nat
\end{lstlisting}\end{quote}
This declares the natural numbers to be the smallest type containing an element, \ttt{zero}, and closed under a successor function, \ttt{succ}. We can then define, say, addition recursively:
\begin{quote}\begin{lstlisting}
definition add : nat → nat → nat
| add x zero := x,
| add x (succ y) := succ (add x y)
\end{lstlisting}\end{quote}
We can also define familiar notation:
\begin{quote}\begin{lstlisting}
notation `ℕ` := nat
notation 0 := zero
notation x `+` y := add x y
\end{lstlisting}\end{quote}
We can then prove that addition is commutative and associative, and go on to define more complex functions and relations, like multiplication, the divisibility relation, and the function which returns the greatest common divisor of two natural numbers. 

To serve that purpose, a proof assistant's language has to support not only writing definitions, but also proofs. Such languages vary widely. In Lean's library, the proof that the greatest common divisor function is commutative runs like this:
\begin{quote}\begin{lstlisting}
theorem gcd.comm (m n : ℕ) : gcd m n = gcd n m :=
dvd.antisymm
  (dvd_gcd !gcd_dvd_right !gcd_dvd_left)
  (dvd_gcd !gcd_dvd_right !gcd_dvd_left) 
\end{lstlisting}\end{quote}
The proof demonstrates the equality of the left- and right-hand side by showing that each one divides the other, using the abstract characterization of \ttt{gcd}. The proof in Isabelle's library does the same thing, but relies on built-in automation to fill in the details:
\begin{quote}\begin{lstlisting}
theorem gcd_comm: "gcd (m::nat) n = gcd n m"
by (auto intro!: dvd.antisym) 
\end{lstlisting}\end{quote}
In contrast, the proof in the SSReflect library is more direct, unfolding the computational definition:
\begin{quote}\begin{lstlisting}
Lemma gcdnC : commutative gcdn.
Proof.
move=> m n; wlog lt_nm: m n / n < m.
  by case: (ltngtP n m) => [||-> //]; last symmetry; auto.
by rewrite gcdnE -{1}(ltn_predK lt_nm) modn_small.
Qed. 
\end{lstlisting}\end{quote}

Proof assistants also have to provide general ``mathematical knowledge management'' and support mathematical conventions. For example, in Lean, the following commands are used to import library developments of the natural numbers and rings, and make information and notation from those libraries readily available:
\begin{quote}\begin{lstlisting}
import data.nat algebra.ring
open nat algebra
\end{lstlisting}\end{quote}
One can globally declare the types of variables used in a particular development:
\begin{quote}\begin{lstlisting}
section
  variables m n k : ℕ
  ...
end
\end{lstlisting}\end{quote}
In ordinary language, this says ``in this section, we use the variables $m$, $n$, and $k$ to range over the natural numbers.'' We can also fix parameters:
\begin{quote}\begin{lstlisting}
section
  parameter G : Group
  variables g₁ g₂ : G
  ...
end
\end{lstlisting}\end{quote}
In words, ``in this section, we fix a group, $G$, and let $g_1$ and $g_2$ range over elements of $G$.''

These examples provide only a small sampling of the ways that users interact with proof assistants. Developing mathematical theories involves not just defining basic objects and proving their properties, but also establishing complex notation and telling the system how to resolve ambiguous expressions, defining algebraic structures and instantiating them, declaring pieces of information to be used by the system's automated procedures, organizing data and structuring complex theories, and tagging facts and data in order to signal intended usage. There is a substantial distance from a foundational definition of addition to the lemmas and theorems of everyday mathematics, and theorem provers are finely engineered to make the passage as smooth as possible.

The point I wish to make here is simply this: whenever someone communicates with a computational proof assistant, they are speaking a formal mathematical language. This language has been designed for a specific purpose, namely, to enable users to develop formal theories smoothly and efficiently. As such, it should be powerful, expressive, and convenient. It should allow users to write mathematical expressions concisely, carry out reasoning steps in an intuitive way, and manage complex information effectively. At the foundational level, the choice of axiomatic system determines the objects that exist and the appropriate rules for reasoning about them, but the higher-level organizational features of the system are no less important: they determine how we can talk about the objects that are guaranteed to exist by the foundational framework, what sorts of conceptual apparatus we can bring to bear upon then, and the styles and patterns of argumentation we can carry out.

\section{The philosophy of mathematics as a design science}
\label{section:philosophy}

In Section~\ref{section:informal}, I described some changes in informal mathematical language that emerged over the course of the nineteenth century, and their effects on mathematical practice. In Section~\ref{section:formal}, I described formal mathematical languages, essentially programming languages, that are used to carry out reasoning in proof assistants, and the role those languages play in structuring our interactions with a system. The parallel should be clear: in both cases we are dealing with the use of language to convey mathematical content and support mathematical reasoning. The languages we use can either serve our purposes well or not, and, in both cases, we continue to tinker with the languages to enable them to serve us better.

The analogy is not perfect, and there are important differences between formal and informal mathematics. Computational proof assistants are developed on the scale of years or decades; mathematical language has been evolving for centuries. The features of a proof assistant are determined by a relatively small design team, responsive, though it may be, to a community of users; changes in mathematical language take place organically across a much larger community. A change to the language of a theorem prover is typically the result of an explicit design decision, and has a precise implementation date; changes to mathematical language are sometimes discussed, but often just happen, sometimes gradually and without fanfare. Perhaps most significantly, the language of a theorem prover has a precisely specified grammar and semantics, and a reference manual that tells you what the various commands do and how they are interpreted by the system. We have to look to more nebulous collections of textbooks and speakers to ascertain the meaning of a mathematical claim. It is a remarkable fact that the mathematical community can communicate with such a high degree of agreement and coherence, without such a reference manual.

Here is the central claim of this essay: when it comes to understanding the power of mathematical language to guide our thought and help us reason well, formal mathematical languages like the ones used by interactive proof assistants provide informative \emph{models} of informal mathematical language. The formal languages underlying foundational frameworks such as set theory and type theory were designed to provide an account of the correct rules of mathematical reasoning, and, as G\"odel observed, they do a remarkably good job. But correctness isn't everything: we want our mathematical languages to enable us to reason efficiently and effectively as well. To that end, we need not just accounts as to what makes a mathematical argument correct, but also accounts of the structural features of our theorizing that help us manage mathematical complexity.

Let me provide some examples of how formal mathematical languages can help. Consider, for example, some of the virtues that I ascribed to the historical restructuring of Dirichlet's proof in Section~\ref{section:informal}: dependencies between components of the proof were minimized, it became easier to understand the arguments and ensure their correctness, and the components themselves became more adaptable, reusable, and generalizable. Computer scientists and programmers will recognize immediately that these are all benefits typically associated with the use of modularity in software design. This is not a coincidence; Section~\ref{section:formal} made the analogies between computer code and formal proof manifest. A formal understanding of modularity, its benefits, and the means by which it can be achieved should therefore illuminate the informal mechanisms by which it is achieved in ordinary mathematics. In a similar way, the means by which programmers manage abstract interfaces that hide the internal details of a program from the outside world should illuminate the informal practices that make mathematical abstraction such a powerful tool.

For another example, consider the act of solving a mathematical problem. This typically requires us to decide, heuristically, which steps from among a combinatorial explosion of options will plausibly move us closer to the goal, and which facts from among a vast store of background knowledge are relevant to the problem at hand. The general problem of using heuristics and contextual cues to navigate an unruly search space is central to artificial intelligence. In informal mathematics, a well-written proof provides subtle yet effective means to manage information and keep the relevant data ready to hand, and a powerful theory is one that provides cues that trigger the steps needed to solve problems in its intended domain. It seems likely that formal modeling and insights from the field of automated reasoning can help us make sense of how informal mathematics manages its problem-solving contexts, and, conversely, informal mathematical data and examples can inform research in AI.

For yet another example, consider the observations above that informal mathematics functions remarkably well without a formal reference manual, and that mathematics has a remarkable way of reinterpreting its past to preserve the validity of prior insights and bodies of knowledge. Here, mathematics has computer science beat: programmers struggle to share algorithms across different code bases, and to preserve code while surrounding infrastructure --- compilers, operating systems, and supporting libraries --- continue to evolve. But once again, the formal mechanisms with which computer scientists cope with the problem can shed light on the mechanisms by which informal mathematics preserves meaning through time and across varied environments, and, conversely, these mechanisms can inform the design of computational systems.

To close this essay, I would like to relate these considerations to the central theme of the essays in this collection, namely, the ontology of mathematics. There are various senses in which ordinary mathematical practice might be said to sanction the existence of an object. Some objects exist only in an imprecise, heuristic sense. When we say that a theorem illuminates an important ``connection'' between algebra and geometry, for example, we do not take ``connections'' to be \emph{bona-fide} mathematical objects with precise properties. We might assign Kepler's points at infinity, described in Jeremy Gray's essay, to this category, if the notion is only used to motivate some of his mathematics without playing a substantial inferential role. 

Next, there are objects that do seem to play a substantial role in our reasoning, but are taken to be \emph{fa\c{c}ons de parler}, convenient shorthands that can be eliminated from ``proper'' mathematical discourse. Some early algebraists seem to have viewed negative magnitudes and nonhomogeneous quantities in this way, so that, for example, an expression like $x^2 + x$ was understood to be the sum of two ``areas,'' $x^2$ and $xy$, where $y$ is a unit length. In 1742, Colin MacLauren wrote \emph{A Treatise of Fluxions}, which aimed to defend the methods of Newton's calculus against Bishop Berkeley's attacks by explaining them in more accessible geometric terms:
\begin{quote}
This is our design in the following treatise; wherein we do not propose to alter Sir {\sc Isaac Newton}'s notion of a fluxion, but to explain and demonstrate his method, by deducing it at length from a few self-evident truths\ldots and, in treating of it, to abstract from all principles and postulates that may require the imagining any other quantities but such as may be easily conceived to have a real existence. 
\end{quote}
The message is this: it is convenient to speak of fluxions, but to be fully rigorous we can proceed in more elementary ways.

At the next level, there are objects that are defined in terms of more basic ones, the way we define a complex number to be pair of real numbers and define a group to be a tuple of data satisfying the group axioms. This, of course, is the most straightforward way to introduce new mathematical objects. 

Finally, there are mathematical objects that exist in the foundational sense that we do not feel the need to define them or explain them in other terms. For the ancients, these included the points, lines, and circles of Euclid's \emph{Elements}; Euclid's ``definitions'' serve to motivate the axioms and diagrammatic rules of inference, but they do not play a direct justificatory role. For most of the history of mathematics, numbers, viewed as discrete or continuous magnitudes, enjoyed such an existence. Today, we generally understand sets and/or functions in this way.

The development of foundational systems like set theory was a watershed in mathematics, providing bedrock foundations that can, in principle, resolve all ontological issues. We no longer need to scratch our heads and wonder whether something exists, since it suffices to produce a set-theoretic definition. And we no longer have to search our intuitions to determine what properties our objects have, since the axioms of set theory serve as the final arbiter.

There are nonetheless good reasons to pursue the broader ontological questions considered in this collection of essays. For one thing, there is still the philosophical question as to what justifies set theory, or any other foundation. Choosing such a foundation amounts to sweeping all one's ontological dust into a corner, but that in and of itself doesn't eliminate it. The best justification for using set theory as an ontology is that we have lots of good reasons to talk about the mathematical objects that we do --- numbers, functions, points on the projective plane, and so on --- and set theory gives us a clean and uniform account of those. In other words, we are justified in talking about sets because we are justified in talking about all those other things, and not the other way around. To round out the story, we need to say something about why we are justified in talking about those other things, and that, in turn, should be based on a meaningful analysis of the inferential and organization role they play in our scientific practices.

Characterizing these roles as ``linguistic'' may convey the impression that mathematical objects are somehow less objective than things we bump into like rocks and chairs, and that we somehow need to explain the ``illusion'' that they are really there. This is nonsense; the fact that we don't bump into something doesn't mean it does not exist. Language works in funny ways, and there are plenty of proper nouns in our language that are not neatly situated in time and space. Any of the buildings of my home university, Carnegie Mellon, may be torn down, and students and faculty come and go, but the university is an entity that somehow survives those changes. We can even imagine the board of trustees shutting down the Pittsburgh operation and relocating to Cleveland, in which case, it would be hard to say exactly what had moved. We can speak about the fifth novel in the Harry Potter series, but we can only point to copies of it. We can produce plenty of red things, but it is hard to point to the color red. How many letters are there in the word ``book''? I can name three, ``b,'' ``o,'' and ``k,'' but some would insist there are four. Perhaps we can disambiguate the question by distinguishing between letters and instances of a letter, but now you have to deal with objects known as ``instances,'' and you still can't point to a letter. And if you manage to bump into the government of the United States, let me know.

It is one thing to ask how we can regiment our informal language to make our claims more precise, for example, for the purposes of scientific inquiry. We may settle on one way of talking about colors in physics and another for the purposes of psychological experimentation, with suitable ways to connect the two. The fact that mathematical objects play a fundamental role across the sciences does not make them any less real than novels, colors, universities, or governments. It seems to me much easier to justify the claim that numbers exist than the claim that novels do.

Nor does aligning our talk of mathematical objects with linguistic norms mean that these norms are ``arbitrary.'' The rules of mathematical discourse are quite rigid, and the language is remarkably stable over time. But, as we have seen, time can even alter the way we talk about fundamental mathematical objects, and it is important to understand the forces that bear upon such changes. To that end, it helps to take a broad view of mathematics, and understand what the objects of mathematics do for us.

Which brings me to an even better reason to take a more expansive view of the ontology of mathematics: mere existence isn't everything. Of all the ``objects'' that set theory or type theory or any other foundation allow us to talk about, only a tiny fraction are of any mathematical interest. So the question is: of all the objects we \emph{can} talk about, which ones \emph{should} we talk about? Trying to answer this question forces us once again to consider the roles that our ontological posits play in the mathematics that we actually do.

In short, it is high time we got over the worry that there is anything inherently mysterious or dubious about mathematical objects, and started paying attention to the things that really matter. On the view I have put forward, even traditional ontological concerns are best addressed by asking why we speak about mathematical objects the way we do, and how the rules and methods we adopt to reason about them provide useful ways of making sense of the world. Questions as to how mathematics enables us to communicate and reason effectively are then of central philosophical importance. The history of mathematics has been an enduring struggle to balance concerns about the meaningfulness, appropriateness, and reliability of our mathematical practices against the drive to adopt innovative methods that enable us to  think better, push our understanding further, and extend our cognitive reach. We need a science that can help us understand the considerations that bear upon our choices of mathematical norms, and the philosophy of mathematics should be that science.

\section*{Notes}

The opening quotation by Wittgenstein is from Part I, \S 165 of his \emph{Remarks on the Foundations of Mathematics}, G.~H.~von Wright, R.~Rhees, G.~E.~M.~Anscombe, eds., MIT Press, 1983. The quotation by Carnap is from his autobiographical chapter in \emph{The Philosophy of Rudolf Carnap}, Paul Schilpp ed., Open Court and Cambridge University Press, 1963. The quotation by Quine is taken from his ``Two Dogmas of Empiricism,'' which first appeared in \emph{The Philosophical Review}, 60: 20--43, 1951 and has been reprinted in many other sources. Alice's exchange with Humpty Dumpty is from Chapter 6 of \emph{Through the Looking Glass}.

The quotation from William Tait is from an essay titled ``The Law of the Excluded Middle and the Axiom of Choice,'' from Alexander George, ed., \emph{Mathematics and Mind}, Oxford University Press, 1994. It is reprinted in William Tait, \emph{The Provenance of Pure Reason: Essays in the Philosophy of Mathematics and Its History}, Oxford University Press, 2005.

The research on the history of Dirichlet's theorem described in Section~\ref{section:informal} can be found in Jeremy Avigad and Rebecca Morris, ``The concept of `character' in Dirichlet’s theorem on primes in an arithmetic progression,'' \emph{Archive for History of Exact Sciences}, 68:265--326, 2014.

G\"odel's incompleteness theorems were first published in ``\"Uber formal unentscheidbare S\"atze der Principia Mathematica und verwandter Systeme, I,'' \emph{Monatshefte f\"ur Mathematik und Physik} 38:173--98, 1931. The excerpt in Section~\ref{section:informal} is taken from the translation in volume I of G\"odel's \emph{Collected Works}, Solomon Feferman et al., eds., Oxford University Press, 1986.

For more on interactive theorem proving and formal verification, see:
\begin{itemize}
 \item Jeremy Avigad and John Harrison, ``Formally verified mathematics,'' \emph{Communications of the ACM}, 57: 66-75, 2014.
 \item Thomas C. Hales, ``Developments in formal proofs,'' \emph{S\'eminaire Bourbaki}, no.~1086, 2014.
 \item Thomas C. Hales, ``Mathematics in the age of the Turing machine,'' in Rod Downey, ed., \emph{Turing's Legacy: Developments from Turing's Ideas in Logic}, ASL Publications, 2014.
 \item The December 2008 special issue on formal proof in the \emph{Notices of the American Mathematical Society}.
\end{itemize}
You can easily find any of the theorem provers mentioned in Section \ref{section:formal} online, where you can browse their libraries and even download the systems and run them on your own computer. Among the snippets of formal text used in Section~\ref{section:formal}, the statement of the prime number theorem is from a formalization I completed in 2004 with the help of Kevin Donnelly, David Gray, and Paul Raff; the statement of the Feit-Thompson theorem is from the ``Mathematical Components'' project, led by Georges Gonthier; the statement of the Kepler conjecture is from the ``Flyspeck'' project, led by Thomas Hales; and the statement of the Blakers-Massey theorem is from a formalization by Kuen-Bang Hou (Favonia), based on a proof by Peter LeFanu Lumsdaine, Eric Finster, and Dan Licata.

Excerpts from MacLauren's \emph{Treatise of Fluxions} can be found in Volume 1 of  \emph{From Kant to Hilbert: A Source Book in the Foundations of Mathematics}, edited by William Ewald, Clarendon Press, 1996.

I may be accused of providing an overly simplified view of ontology. My discussion of the difficulties in making sense of everyday natural language borrows examples from Gilbert Ryle's \emph{The Concept of Mind} (University of Chicago Press, 1949), W.~V.~O.~Quine's \emph{Word and Object} (MIT Press, 1960), and John Burgess and Gideon Rosen's \emph{A Subject with No Object: Strategies for Nominalistic Interpretation of Mathematics} (Oxford University Press, 1999). All are well worth reading.

I am grateful to Ernest Davis, Jeremy Gray, and Robert Lewis for discussions and helpful comments, and to Eric Tressler for corrections.

\end{document}